\documentclass{amsart}
\usepackage{amssymb,amscd}

\newtheorem{lemma}{Lemma}[section]
\newtheorem{theorem}[lemma]{Theorem}
\newtheorem{proposition}[lemma]{Proposition}

\newtheorem{problem}[lemma]{Problem}

\theoremstyle{definition}

\newcommand{\F}{\mathcal {F}}
\newcommand{\Ra}{\Rightarrow}
\newcommand{\U}{\mathcal U}
\newcommand{\V}{\mathcal V}

\newcommand{\A}{\mathcal A}
\newcommand{\Z}{\mathcal Z}

\newcommand{\IN}{\mathbb N}
\newcommand{\IZ}{\mathbb Z}
\newcommand{\IC}{\mathbb C}

\newcommand{\invG}{\overset{\leftrightarrow}{G}}

\newcommand{\inv[1]}{\overset{\leftrightarrow}{#1}}

\title[Algebra in superextensions of groups, I]{Algebra in superextensions of groups, I:\\ zeros and commutativity}
\author{T.Banakh, V.Gavrylkiv, O.Nykyforchyn}
\address{Ivan Franko National University of Lviv, Ukraine}
\email{tbanakh@yahoo.com}
\address{Vasyl Stefanyk Precarpathian National University, Ivano-Frankivsk, Ukraine}
\email{vgavrylkiv@yahoo.com}

\begin{document}
\begin{abstract}Given a group $X$ we study the algebraic structure of its superextension
 $\lambda(X)$. This is a right-topological semigroup consisting of all maximal linked systems
 on $X$ endowed with the operation $$\A\circ\mathcal B=\{C\subset X:\{x\in X:x^{-1}C\in\mathcal B\}\in\A\}$$
 that extends the group operation of $X$.
We characterize right zeros of $\lambda(X)$ as invariant maximal
linked systems on $X$ and prove that
 $\lambda(X)$ has a right zero if and only if each element of $X$ has odd order. On the other hand,
  the semigroup $\lambda(X)$ contains a left zero if and only if it contains a zero if and only if
   $X$ has odd order $|X|\le5$. The semigroup $\lambda(X)$ is commutative if and only if $|X|\le4$.
   We finish the paper with a complete description of the algebraic structure of the semigroups $\lambda(X)$
  for all groups $X$ of cardinality $|X|\le5$.
\end{abstract}

\maketitle
\setcounter{tocdepth}{1}
\tableofcontents

\section*{Introduction}

 After the topological proof of the Hindman theorem \cite{Hind} given by Galvin and Glazer\footnote{Unpublished, see
\cite[p.102]{HS}, \cite{H2}},
topological methods become a standard tool in the modern
combinatorics of numbers, see \cite{HS}, \cite{P}. The crucial
point is that any semigroup operation $\ast$ defined on a
discrete space $X$ can be extended to a right-topological
semigroup operation on $\beta(X)$, the Stone-\v Cech
compactification of $X$. The extension of the operation from $X$
to $\beta(X)$ can be defined by the simple formula:
\begin{equation}\label{extension}\U\circ\V=\{A\subset X:\{x\in X:x^{-1}A\in\V\}\in\U\},
\end{equation}
where $\U,\V$ are ultrafilters on $X$ and $x^{-1}A=\{y\in X:xy\in A\}$. Endowed with the so-extended
operation, the Stone-\v Cech compactification $\beta(X)$ becomes a
compact right-topological semigroup. The algebraic properties of
this semigroup (for example, the existence of idempotents or
minimal left ideals) have important consequences in combinatorics
of numbers, see \cite{HS}, \cite{P}.

The Stone-\v Cech compactification $\beta(X)$ of $X$ is the
subspace of the double power-set $\mathcal P(\mathcal P(X))$,
which is a complete lattice with respect to the operations of
union and intersection. In \cite{G2} it was observed that the
semigroup operation extends not only to $\beta(X)$ but also to the
complete sublattice $G(X)$ of $\mathcal P(\mathcal P(X))$
generated by $\beta(X)$. This complete sublattice consists of all
inclusion hyperspaces over $X$.

By definition, a family $\F$ of non-empty subsets of a discrete
space $X$ is called an {\em inclusion hyperspace} if $\F$ is
monotone in the sense that a subset $A\subset X$ belongs to $\F$
provided $A$ contains some set $B\in\F$. On the set $G(X)$ there
is an important transversality operation assigning to each
inclusion hyperspace $\F\in G(X)$ the inclusion hyperspace
$$\F^\perp=\{A\subset X:\forall F\in\F\;(A\cap F\ne\emptyset)\}.$$
This operation is involutive in the sense that $(\F^\perp)^\perp=\F$.

It is known that the family $G(X)$ of inclusion hyperspaces on $X$
is closed in the double power-set $\mathcal P(\mathcal
P(X))=\{0,1\}^{\mathcal P(X)}$ endowed with the natural product
topology.  The induced topology on $G(X)$ can be described
directly: it is generated by the sub-base consisting of the sets
$$U^+=\{\F\in G(X):U\in\F\}\mbox{ and }U^-=\{\F\in G(X):U\in\F^\perp\}$$
where $U$ runs over subsets of $X$. Endowed with this topology, $G(X)$
becomes a Hausdorff supercompact  space. The latter means that each cover
 of $G(X)$ by the sub-basic sets has a 2-element subcover.

The extension of a binary operation $\ast$ from $X$ to $G(X)$ can
be defined in the same way as for ultrafilters, i.e., by the
formula~(\ref{extension}) applied to any two inclusion hyperspaces
$\U,\V\in G(X)$. In \cite{G2} it was shown that for an associative
binary operation $\ast$ on $X$ the space $G(X)$ endowed with the
extended operation becomes a compact right-topological semigroup. The algebraic
 properties of this semigroups were studied in details in \cite{G2}.

Besides the Stone-\v Cech compactification $\beta(X)$, the semigroup $G(X)$ contains
many important spaces as closed subsemigroups. In particular,  the
space
$$\lambda(X)=\{\F\in G(X):\F=\F^\perp\}$$of maximal linked systems on
$X$ is a closed subsemigroup of $G(X)$. The space $\lambda(X)$ is
well-known in  General and Categorial Topology as the {\em superextension}
of $X$, see \cite{vM}, \cite{TZ}.  Endowed with the extended binary operation, the
superextension $\lambda(X)$ of a semigroup $X$ is a supercompact
right-topological semigroup containing $\beta(X)$ as a subsemigroup.

The space $\lambda(X)$ consists of maximal linked systems on $X$. We recall that a system of
subsets $\mathcal L$ of $X$ is {\em linked} if $A\cap B\ne\emptyset$ for all $A,B\in\mathcal L$.
An inclusion hyperspace $\A\in G(X)$ is linked if and only if $\A\subset\A^\perp$. The family of
all linked inclusion hyperspace on $X$ is denoted by $N_2(X)$. It is a closed subset in $G(X)$.
Moreover, if $X$ is a semigroup, then $N_2(X)$ is a closed subsemigroup of $G(X)$. The superextension
 $\lambda(X)$ consists of all maximal elements of $N_2(X)$, see \cite{G1}, \cite{G2}.

In this paper we start a systematic investigation of the algebraic structure of the semigroup $\lambda(X)$.
This program will be continued in the forthcoming papers \cite{BG2} and \cite{BG3}.
The interest to studying the semigroup $\lambda(X)$ was motivated by the fact that for each maximal linked
system $\mathcal L$ on $X$ and each partition $X=A\cup B$ of $X$ into two sets $A,B$ either $A$ or $B$
belongs to $\mathcal L$. This makes possible to apply maximal linked systems to Combinatorics and Ramsey Theory.

In this paper we concentrate on describing zeros and commutativity of the semigroup $\lambda(X)$.
In Proposition~\ref{rightzero} we shall show that a maximal linked system $\mathcal L\in\lambda(X)$ is
 a right zero of $\lambda(X)$  if and only if $\mathcal L$ is invariant in the sense that $xL\in\mathcal L$
 for all $L\in\mathcal L$ and all $x\in X$. In Theorem~\ref{odd} we shall prove that a group $X$ admits an
  invariant maximal linked system (equivalently, $\lambda(X)$ contains a right zero) if and only if each
  element of $X$ has odd order. The situation with (left) zeros is a bit different: a maximal linked
  system $\mathcal L\in \lambda(X)$ is a left zero in $\lambda(X)$ if and only if $\mathcal L$ is a zero
  in $\lambda(X)$ if and only if $\mathcal L$ is a unique invariant maximal linked system on $X$.
   The semigroup $\lambda(X)$ has a (left) zero if and only if $X$ is a finite group of odd order $|X|\le 5$
   (equivalently, $X$ is isomorphic to the cyclic group $C_1$, $C_3$ or $C_5$). The semigroup $\lambda(X)$
   rarely is commutative: this holds if and only if the group $X$ has finite order $|X|\le 4$.

We start the paper studying self-linked subsets of groups. By
definition, a subset $A$ of a group $X$ is called {\em
self-linked} if $A\cap xA\ne\emptyset$ for all $x\in X$.
 In Proposition~\ref{slbound} we shall give lower and upper bounds for the smallest cardinality $sl(X)$
 of a self-linked subset of $X$. We use those bounds to characterize groups $X$ with $sl(X)\ge |X|/2$
 in Theorem~\ref{upbound}.

In Section~\ref{maxinv} we apply self-linked sets to  evaluating
the cardinality of the (rectangular) semigroup $\inv[\lambda](X)$
of maximal invariant linked systems on a group $X$. In
Theorem~\ref{t2.2}
 we show that for an infinite group $X$ the cardinality of $\inv[\lambda](X)$ equals $2^{2^{|X|}}$.
 In Proposition~\ref{p2.3} and
Theorem~\ref{fingr} we calculate the cardinality of
$\inv[\lambda](X)$ for all finite groups $X$ of
 order $|X|\le 8$ and also detect groups $X$ with $|\inv[\lambda](X)|=1$. In Sections~\ref{lzeros}
 and \ref{scom} these results are applied for characterizing groups $X$ whose superextensions have
 zeros or are commutative.

We finish the paper with a description of the algebraic structure
of the superextensions of groups $X$ of order $|X|\le 5$.

Now a couple of words about notations. Following the algebraic
tradition, by $C_n$ we denote the cyclic group of order $n$ and by
$D_{2n}$ the dihedral group of cardinality $2n$, that is, the
isometry group of the regular $n$-gon. For a group $X$ by $e$ we
denote the neutral element of $X$. For a real number
 $x$ we put $$\lceil x\rceil=\min\{n\in\IZ:n\ge x\}\mbox{ and }\lfloor x\rfloor=\max\{n\in\IZ:n\le x\}.$$

\section{Self-linked sets in groups}

In this section we study self-linked subsets in groups. By definition, a subset $A$ of a group $G$ is
 {\em self-linked} if $A\cap xA\ne\emptyset$ for each $x\in G$. In fact, this notion can be defined in
 the more general context of $G$-spaces.

By a {\em $G$-space} we understand a set $X$ endowed with a left action $G\times X\to X$ of a group $G$.
Each group $G$ will be considered as a $G$-space endowed with the left action of $G$. An important example
of a $G$-space is the homogeneous space $G/H=\{xH:x\in G\}$ of a group $G$ by a subgroup $H\subset G$.

A subset $A\subset X$ of a $G$-space $X$ defined to be {\em self-linked} if $A\cap gA\ne\emptyset$ for
 all $g\in G$.  Let us observe that a subset $A\subset G$ of a group $G$ is self-linked if and only if $AA^{-1}=G$.

For a $G$-space $X$ by $sl(X)$ we denote the smallest cardinality $|A|$ of a self-linked subset $A\subset X$.
Some lower and upper bounds for $sl(G)$ are established in the following proposition.

\begin{proposition}\label{slbound} Let $G$ be a finite group and $H$ be a subgroup of $G$. Then
\begin{enumerate}
\item $sl(G)\ge (1+\sqrt{4|G|-3})/2$;
\item $sl(G)\le sl(H)\cdot sl(G/H)\le sl(H)\cdot\lceil(|G/H|+1)/2\rceil.$
\item $sl(G)<|H|+|G/H|$.
\end{enumerate}
\end{proposition}

\begin{proof}
1. Take any self-linked set $A\subset G$ of cardinality $|A|=sl(G)$ and consider the surjective
map $f:A\times A\to G$, $f:(x,y)\mapsto xy^{-1}$. Since $f(x,y)=xy^{-1}=e$ for all
$(x,y)\in\Delta_A=\{(x,y)\in A^2:x=y\}$, we get $|G|=|G\setminus\{e\}|+1\le |A^2\setminus\Delta_A|+1=sl(G)^2-sl(G)+1$,
which just implies that $sl(G)\ge (1+\sqrt{4|G|-3})/2$.
\smallskip

2a. Let $H$ be a subgroup of $G$. Take self-linked sets $A\subset H$ and $\mathcal B\subset G/H=\{xH:x\in G\}$
 having sizes $|A|=sl(H)$ and $|\mathcal B|=sl(G/H)$. Fix any subset $B\subset G$ such that $|B|=|\mathcal B|$
 and $\{xH:x\in B\}=\mathcal B$.
We claim that the set $C=BA$ is self-linked. Given arbitrary $x\in G$ we should prove that the intersection
$C\cap xC$ is not empty. Since $\mathcal B$ is self-linked, the intersection $\mathcal B\cap x\mathcal B$
contains the coset $bH=xb'H$ for some $b,b'\in B$. It follows that $b^{-1}xb'\in H=AA^{-1}$. The latter equality
follows from the fact that the set $A\subset H$ is self-linked in $H$. Consequently, $b^{-1}xb'=a'a^{-1}$ for
some $a,a'\in A$. Then $xC\ni xb'a=ba'\in C$ and thus $C\cap xC\ne\emptyset$. The self-linkedness of $C$  implies
the desired upper bound$$sl(G)\le|C|\le|A|\cdot|B|=sl(H)\cdot sl(G/H).$$
\smallskip

2b. Next, we show that $sl(G/H)\le\lceil(|G/H|+1)/2\rceil$. Take any subset $A\subset G/H$ of size
$|A|=\lceil(|G/H|+1)/2\rceil$ and note that $|A|>|G/H|/2$. Then for each $x\in G$ the shift $xA$ has
size $|xA|=|A|>|G/H|/2$. Since $|A|+|xA|>|G/H|$, the sets $A$ and $xA$ meet each other. Consequently, $A$ is
self-linked and $sl(G/H)\le|A|=\lceil(|G/H|+1)/2\rceil$.
\smallskip

3. Pick a subset $B\subset G$ of size $|B|=|G/H|$ such that $BH=G$ and observe that the set $A=H\cup B$ is
self-linked and has size $|A|\le|H|+|B|-1$ (because $B\cap H$ is a singleton).
\end{proof}

\begin{theorem}\label{upbound}
For a finite group $G$
\begin{enumerate}
\item[(i)] $sl(G)=\lceil(|G|+1)/2\rceil>|G|/2$ if and only if $G$ is isomorphic to one of the groups: $C_1$, $C_2$, $C_3$, $C_4$, $C_2\times C_2$, $C_5$, $D_6$, $(C_2)^3$;
\item[(ii)] $sl(G)=|G|/2$ if and only if $G$ is isomorphic to one of the groups: $C_6$, $C_8$, $C_4\times C_2$, $D_8$, $Q_8$.
\end{enumerate}
\end{theorem}

\begin{proof} I. First we establish the inequality $sl(G)<|G|/2$ for all groups $G$  not isomorphic to
the groups appearing in the items (i), (ii). Given such a group $G$ we should find a self-linked subset
$A\subset G$ with $|A|<|G|/2$.

We consider 8 cases.

1) $G$ contains a subgroup $H$ of order $|H|=3$ and index $|G/H|=3$. Then $sl(H)=2$ and we can apply
Proposition~\ref{slbound}(2) to conclude that
$$sl(G)\le sl(H)\cdot sl(G/H)\le 2\cdot 2<9/2=|G|/2.$$

2) $|G|\notin\{9,12,15\}$ and $G$ contains a subgroup $H$ of order $n=|H|\ge 3$ and index $m=|G/H|\ge 3$.
It this case $n+m-1<nm/2$ and $sl(G)\le|H|+|G/H|-1=n+m-1<nm/2$ by Proposition~\ref{slbound}(3).

3) $G$ is cyclic of order $n=|G|\ge 9$. Given a generator $a$ of
$G$, construct a sequence $(x_i)_{2\le i\le n/2}$ letting
$x_2=a^0$, $x_3=a$, $x_4=a^3$, $x_5=a^5$, and $x_{i}=x_{i-1}a^{i}$
for $5<i\le n/2$. Then the set $A=\{x_i:2\le i\le n/2\}$ has size
$|A|<n/2$ and is self-linked.

4) $G$ is cyclic of order $|G|=7$. Given a generator $a$ of $G$ observe that $A=\{e,a,a^3\}$ is a 3-element
self-linked subset and thus $sl(G)\le 3<|G|/2$.

5) $G$ contains a cyclic subgroup $H\subset G$ of prime order $|H|\ge 7$. By the preceding two cases, $sl(H)<|H|/2$
and then $sl(G)\le sl(H)\cdot sl(G/H)< \frac{|H|}2\cdot\frac{|G|}{|H|}=|G|/2.$

6) $|G|>6$ and $|G|\notin\{8,10,12\}$. If $|G|$ is prime or $|G|=15$, then $G$ is cyclic of order $|G|\ge 7$ and
thus has $sl(G)<|G|/2$ by the items (3), (4). If $|G|=2p$ for some prime number $p$, then $G$ contains a cyclic
 subgroup of order $p\ge 7$ and thus has $sl(G)<|G|/2$ by the item (5). If $|G|=4n$ for some $n\ge 4$, then by Sylow's
 Theorem (see \cite[p.74]{OA}), $G$ contains a subgroup $H\subset G$ of order $|H|=4$ and index $|G/H|\ge 4$.
 Then $sl(G)<|G|/2$ by the item (2).
If the above conditions do not hold, then $|G|=nm\ne 15$ for some odd numbers $n,m\ge3$ and we can apply the items
 (1) and (2) to conclude that $sl(G)<|G|/2$.

7) If $|G|=8$, then $G$ is isomorphic to one of the groups: $C_8$, $C_2\times C_4$, $(C_2)^3$, $D_8$, $Q_8$. All those groups appear in the items (i), (ii) and thus are excluded from our consideration.

8) If $|G|=10$, then $G$ is isomorphic to $C_{10}$ or $D_{10}$. If $G$ is isomorphic to $C_{10}$, then
$sl(G)<|G|/2$ by the item (3). If $G$ is isomorphic to $D_{10}$, then $G$ contains an element $a$ of order 5
and an element $b$ of order 2 such that  $bab^{-1}=a^{-1}$. Now it is easy to check that the 4-element set
$A=\{e,a,b,ba^2\}$ is self-linked and hence $sl(G)\le 4<|G|/2$.

9) In this item we consider groups $G$ with $|G|=12$. It is well-known that there are five non-isomorphic
groups of order 12: the cyclic group $C_{12}$, the direct sum of two cyclic groups $C_6\oplus C_2$, the
dihedral group $D_{12}$, the alternating group $A_4$, and the semidirect product $C_3\rtimes C_4$ with
 presentation $\langle a,b|a^4=b^3=1, aba^{-1}=b^{-1}\rangle$.

If $G$ is isomorphic to $C_{12}$, $C_6\oplus C_2$ or $A_4$, then $G$ contains a normal 4-element subgroup $H$.
 By Sylow's Theorem, $G$ contains also an element $a$ of order 3. Taking into account that $a^2\notin H$ and
 $Ha^{-1}=a^{-1}H$, we conclude that the 5-element set $A=\{a\}\cup H$ is self-linked and hence $sl(G)\le 5<|G|/2$.

If $G$ is isomorphic to $C_3\rtimes C_4$, then $G$ contains a normal subgroup $H$ of order 3 and an element
 $a\in G$ such that $a^2\notin H$. Observe that the 5-element set $A=H\cup\{a,a^2\}$ is self-linked. Indeed,
 $AA^{-1}\supset H\cup aH\cup a^2H\cup Ha^{-1}=G$. Consequently, $sl(G)\le 5<|G|/2$.

Finally, consider the case of the dihedral group $D_{12}$. It contains an element $a$ generating a cyclic subgroup
 of order 6 and an element $b$ of order 2 such that $bab^{-1}=a^{-1}$. Consider the 5-element set $A=\{e,a,a^3,b,ba\}$
  and note that
$AA^{-1}=\{e,a,a^3,b,ba\}\cdot\{e,a^5,a^3,b,ba\}=G$. This yields the desired inequality $sl(G)\le 5<6=|G|/2$.
\smallskip

Therefore we have completed the proof of the inequality $sl(G)<|G|/2$ for all groups not appearing in the items
(i),(ii) of the theorem.
\smallskip

II. Now we shall prove the item (i).

The lower bound from Proposition~\ref{slbound}(1) implies that $sl(G)=\lceil(|G|+1)/2\rceil>|G|/2$ for all
groups $G$ with $|G|\le 5$.

It remains to check that $sl(G)>|G|/2$ if $G$ is isomorphic to
$D_6$ or $C_2^3$. First we consider the case $G=D_6$. In this case
$G$ contains a normal 3-element subgroup $T$. Assuming that
$sl(G)\le |G|/2=3$, find a self-linked 3-element subset $A$.
Without loss of generality we can assume that the neutral element
$e$ of $G$ belongs to $A$ (otherwise replace $A$ by a suitable
shift $xA$).
  Taking into account that  $AA^{-1}=G$, we conclude that $A\not\subset T$ and thus we can
   find an element $a\in A\setminus T$. This element has order 2. Then
$$AA^{-1}=\{e,a,b\}\cdot\{e,a,b^{-1}\}=\{e,a,b,a,e,ba,b^{-1},ba,e\}\ne G,$$which is a contradiction.

Now assume that $G$ is isomorphic to $C_2^3$.  In this case $G$ is
the 3-dimensional linear space over the field $C_2$. Assuming that
$sl(A)\le 4=|G|/2$, find a 4-element self-linked subset $A\subset
G$. Replacing $A$ by a suitable shift, we can assume that  $A$
contains a neutral element $e$ of $G$. Since $AA^{-1}=G$, the set
$A$ contains three linearly independent points $a,b,c$. Then
 $$AA^{-1}=\{e,a,b,c\}\cdot\{e,a,b,c\}=\{e,a,b,c,ab,ac,bc\}\ne G,$$which contradicts the choice of $A$.
\smallskip

III. Finally, we prove the equality $sl(G)=|G|/2$ for the groups appearing in the item (ii).

If $G=C_6$, then $sl(G)\ge 3$ by Proposition~\ref{slbound}(1). On the other hand, we can check that
 for any generator $a$ of $G$ the 3-element subset $A=\{e,a,a^3\}$ is self-linked in $G$, which yields $sl(G)=3=|G|/2$.

If $|G|=8$, then $sl(G)\ge 4$ by Proposition~\ref{slbound}(1).

If $G$ is cyclic of order 8 and $a$ is a generator of $G$, then the set $A=\{e,a,a^3,a^4\}$ is
self-linked and thus $sl(C_8)=4$.

If $G$ is isomorphic to $C_4\oplus C_2$, then $G$ has two commuting generators $a,b$ such that $a^4=b^2=1$.
One can check that the set
$A=\{e,a,a^2,b\}$ is self-linked and thus $sl(C_4\oplus C_2)=4$.

If $G$ is isomorphic to the dihedral group $D_8$, then $G$ has two generators $a,b$ connected by the relations
$a^4=b^2=1$ and $bab^{-1}=a^{-1}$. One can check that the 4-element subset $A=\{e,a,b,ba^2\}$ is self-linked.

If $G$ is isomorphic to the group $Q_8=\{\pm1,\pm i,\pm j,\pm k\}$ of quaternion units, then we can check that
the 4-element subset $A=\{-1,1,i,j\}$ is self-linked and thus $sl(Q_8)=4$.
\end{proof}

In the following proposition we complete Theorem~\ref{upbound} calculating the values of the cardinal
 $sl(G)$ for all groups $G$ of cardinality $|G|\le 13$.

\begin{proposition}\label{fincar1} The number $sl(G)$ for a group $G$ of size
$|G|\le 13$ can be found from the table:
{\small
$$
\begin{array}{|c|c|c|c|cc|cc|ccccc|}
\hline
G& C_2& C_3 & C_5&C_4 &C_2\oplus C_2 &C_6& D_6 & C_8 &C_2\oplus C_4  & D_8 & Q_8&C_2^3\\
\hline
sl(G)&2&2&3&3&3&3&4&4&4&4&4&5\\
\hline
\hline
G&C_{7}& C_{11}& C_{13} &C_9&C_3\oplus C_3 &C_{10}& D_{10} & C_{12} &C_2\oplus C_6 &D_{12}& A_4 & C_3\rtimes C_4\\
\hline
sl(G)&3&4&4&4&4&4&4&4&5&5&5&5\\
\hline
\end{array}$$
}
\end{proposition}

\begin{proof} For groups $G$ of order $|G|\le 10$ the value of $sl(G)$ is uniquely determined by the
 lower bound $sl(G)\ge\frac{1+\sqrt{4|G|-3}}2$ from Proposition~\ref{slbound}(1) and the upper bound
 from Theorem~\ref{upbound}. It remains to consider the groups $G$ of order $11\le |G|\le 13$.
\smallskip

1. If $G$ is cyclic of order 11 or 13, then take a generator $a$ of $G$ and check that the 4-element
 set $A=\{e,a^4,a^5,a^7\}$ is self-linked, witnessing that $sl(C_{12})=4$.
\smallskip

2. If $G$ is cyclic of order 12, then take a generator $a$ for $G$ and check that the 4-element subset
$A=\{e,a,a^3,a^7\}$ is self-linked witnessing that $sl(G)=4$.
\smallskip

It remains to consider all other groups of order 12. Theorem~\ref{upbound} gives us an upper bound
$sl(G)\le 5$. So, we need to show that $sl(G)>4$ for all non-cyclic groups $G$ with $|G|=12$.
\smallskip

3. If $G$ is isomorphic to $C_6\oplus C_2$ or $A_4$, then $G$
contains a normal subgroup $H$ isomorphic to $C_2\oplus C_2$.
Assuming that $sl(G)=4$, we can find a 4-element self-linked
subset $A\subset G$. Since $AA^{-1}=G$, we can find a suitable
shift $xA$ such that  $xA\cap H$ contains the neutral element $e$
of $G$ and some other element $a$ of $H$. Replacing $A$ by $xA$,
we can assume that $e,a\in A$. Since $A\not\subset H$, there is a
point $b\in A\setminus H$. Since the quotient group $G/H$ has
order 3, $bH\cap Hb^{-1}=\emptyset$.

Concerning the forth element $c\in A\setminus\{e,a,b\}$ there are three possibilities: $c\in H$, $c\in b^{-1}H$,
and $c\in bH$.
If $c\in H$, then $bH=bH\cap AA^{-1}=b(A\cap H)^{-1}$ consists of 3 elements which is a contradiction.
If $c\in b^{-1}H$, then $H=H\cap AA^{-1}=\{e,a\}$, which is absurd.
So, $c\in bH$ and thus $c=bh$ for some $h\in H$. Since $h=h^{-1}$, we get  $cb^{-1}=bhb^{-1}=bh^{-1}b^{-1}=bc^{-1}$.
Then $H=H\cap AA^{-1}=\{e,a,cb^{-1},bc^{-1}\}$ has cardinality $|H|=|\{e,a,cb^{-1}=bc^{-1}\}|\le 3$, which is not true.
This contradiction completes the proof of the inequality $sl(G)>4$ for the groups $C_6\oplus C_2$ and $A_4$.
\smallskip

4. Assume that $G$ is isomorphic to the dihedral group $D_{12}$. Then $G$ contains a normal cyclic subgroup $H$
of order 6, and for each $b\in G\setminus H$ and $a\in H$ we get $b^2=e$ and $bab^{-1}=a^{-1}$. Assuming
that $sl(D_{12})=4$, we can find a 4-element self-linked subset $A\subset G$. Let $a$ be a generator of the
group $H$. Since $a\in AA^{-1}=G$, we can find two element $x,y\in A$ such that $a=xy^{-1}$. Then the shift
$Ay^{-1}$ contains $e$ and $a$. Replacing $A$ by $Ay^{-1}$, if necessary, we can assume that $e,a\in A$.
Since $A\not\subset H$, there is an element $b\in A\setminus H$. Concerning the forth element
$c\in A\setminus\{e,a,b\}$ there are two possibilities: $c\in H$ and $c\notin H$. If $c\in H$, then the set
 $A_H=A\cap H=\{e,a,c\}$ contains three elements and is equal to $bA_H^{-1}b^{-1}$, which implies
 $bA_H^{-1}=A_Hb^{-1}=bA_H^{-1}\cup A_Hb^{-1}=AA^{-1}\cap bH=bH$. This is a contradiction, because
 $|H|=4>3=|bA_H^{-1}|$. Then $c\in bH$ and hence $H=H\cap AA^{-1}=\{e,a,a^{-1},bc^{-1},cb^{-1}\}$ which
 is not true because $|H|=6>5$.
\smallskip

5. Assume that $G$ is isomorphic to the semidirect product $C_3\rtimes C_4$ and hence has a presentation
$\langle a,b|a^4=b^3=1, aba^{-1}=b^{-1}\rangle$. Then the cyclic subgroup $H$ generated by $b$ is normal
in $G$ and the quotient $G/H$ is cyclic of order 4. Assuming that $sl(G)=4$, take any 4-element self-linked
subset $A\subset G$.

After a suitable shift of $A$, we can assume that $e,b\in A$. Since $A\not\subset H$, there is an element
$c\in A\setminus H$. We claim that the fourth element
$d\in A\setminus\{e,b,c\}$ does not belong to $H\cup cH\cup c^{-1}H$. Otherwise,
$AA^{-1}\subset H\cup cH\cup c^{-1}H\ne G$. This implies that one of the elements, say $c$ belongs to the
coset $a^2H$ and the other to $aH$ or $a^{-1}H$. We lose no generality assuming that $d\in aH$.
Then $c=a^2b^i$, $d=ab^j$ for some $i,j\in\{-1,0,1\}$.
It follows that $$aH=aH\cap AA^{-1}=
\{d, db^{-1},cd^{-1}\}=
\{ab^j,ab^{j-1},a^2b^{i-j}a^{-1}\}
=\{ab^j, ab^{j-1},ab^{j-i}\}$$which implies that $i=-1$ and thus $c=a^2b^{-1}$.
In this case we arrive to a contradiction looking at
$$a^2H\cap AA^{-1}=\{c,cb^{-1},c^{-1},bc^{-1}\}=\{a^2b^{-1},a^2b^{-2},ba^2,b^2a^2\}\not\ni a^2.$$
\end{proof}

\begin{problem} What is the value of $sl(G)$ for other groups $G$ of small cardinality?
Is $sl(G)=\lceil(1+\sqrt{4|G|-3})/2\rceil$ for all finite cyclic groups $G$?
\end{problem}

\section{Maximal invariant linked systems}\label{maxinv}

In this section we study (maximal) invariant linked systems on groups. An
inclusion hyperspace
 $\A$ on a group $X$ is called {\em invariant} if $x\A=\A$ for all
  $x\in X$. The set of all invariant inclusion hyperspaces on $X$ is denoted by $\invG(X)$.
  By \cite{G2}, $\inv[G](X)$ is a closed rectangular subsemigroup of $G(X)$ coinciding with the minimal
  ideal of $G(X)$. The {\em rectangularity} of $\inv[G](X)$ means that $\A\circ\mathcal B=\mathcal B$
  for all $\A,\mathcal B\in\inv[G](X)$.

  Let $\inv[N]_2(X)=N_2(X)\cap\inv[G](X)$ denote the set of all invariant linked systems
  on $X$ and $\inv[\lambda](X)=\max\inv[N]_2(X)$ be the family of all maximal elements of
  $\inv[N]_2(X)$. Elements of $\inv[\lambda](X)$ are called {\em maximal invariant linked systems}.
  The reader should be concisions of the fact that maximal invariant linked systems need not be maximal linked!

\begin{theorem}\label{t2.1} For every group $X$ the set $\inv[\lambda](X)$ is a non-empty closed rectangular
subsemigroup of $G(X)$.
\end{theorem}

\begin{proof} The rectangularity of $\inv[\lambda](X)$ implies from the rectangularity of $\inv[G](X)$
 established in \cite[\S5]{G2} and the inclusion $\inv[\lambda](X)\subset\inv[G](X)$.

The Zorn Lemma implies that each invariant linked system on $X$
(in particular, $\{X\}$) can be enlarged to a maximal invariant linked system on $X$.
 This observation implies the set $\inv[\lambda](X)$ is not empty. Next, we show that
the subsemigroup  $\inv[\lambda](X)$ is closed in $G(X)$.
Since the set $\inv[N]_2(X)=N_2(X)\cap\inv[G](X)$ is closed in $G(X)$, it suffices to show that
$\inv[\lambda](X)$ is closed in $\inv[N]_2(X)$. Take any invariant linked system
$\mathcal L\in\inv[N]_2(X)\setminus\inv[\lambda](X)$. Being not maximal
invariant, the linked system $\mathcal L$ can be enlarged to a maximal invariant linked system
$\mathcal M$ that contains a subset $B\in\mathcal M\setminus\mathcal L$. Since $\mathcal M\ni B$
 is invariant, the system $\{xB:x\in X\}\subset\mathcal M$ is linked. Observe that $B\notin\mathcal L$
 and $B\in\mathcal M\supset\mathcal L$ implies $X\setminus B\in\mathcal L^\perp$ and $B\in\mathcal L^\perp$.
 We claim that $O(\mathcal L)=B^-\cap(X\setminus B)^-\cap\inv[N]_2(X)$ is a neighborhood of
  $\mathcal L$ in $\inv[N]_2(X)$ that misses the set $\inv[\lambda](X)$. Indeed, for any
  $\A\in O(\mathcal L)$, we get that $\A$ is an invariant linked system such that $B\in\A^\perp$.
  Observe that for every $x\in X$ and $A\in\A$ we get $x^{-1}A\in\A$ by the invariantness of
   $\A$ and hence the set $B\cap x^{-1}A$ and its shift $xB\cap A$ both are not empty.
    This witnesses that $xB\in\A^\perp$ for every $x\in X$. Then the maximal invariant linked system
    generated by $ \A\cup\{xB:x\in X\}$ is an invariant linked enlargement of $\A$,
    which shows that $\A$ is not maximal invariant linked.
\end{proof}

Next, we shall evaluate the cardinality of $\inv[\lambda](X)$.

\begin{theorem}\label{t2.2} For any infinite group $X$ the semigroup $\inv[\lambda](X)$ has cardinality
 $|\inv[\lambda](X)|=2^{2^{|X|}}$.
\end{theorem}

\begin{proof} The upper bound $|\inv[\lambda](X)|\le2^{2^{|X|}}$ follows from the chain
of inclusions: $$\inv[\lambda](X)\subset G(X)\subset \mathcal P(\mathcal P(X)).$$

Now we prove that $|\inv[\lambda](X)|\ge 2^{2^{|X|}}$. Let $|X|=\kappa$ and $X=\{x_\alpha:\alpha<\kappa\}$
be an injective enumeration of $X$ by ordinals $<\kappa$ such that $x_0$ is the neutral element of $X$.
 For every $\alpha<\kappa$ let $B_\alpha=\{x_\beta,x_\beta^{-1}:\beta<\alpha\}$. By transfinite induction,
 choose a transfinite sequence $(a_\alpha)_{\alpha<\kappa}$ such that $a_0=x_0$ and
$$a_\alpha\notin B_\alpha^{-1}B_\alpha A_{<\alpha}$$where $A_{<\alpha}=\{a_\beta:\beta<\alpha\}$.

Consider the set $A=\{a_\alpha:\alpha<\kappa\}$. By \cite[3.58]{HS}, the set $U_\kappa(A)$ of $\kappa$-uniform
 ultrafilters on $A$ has cardinality $|U_\kappa(A)|=2^{2^{\kappa}}$. We recall that an ultrafilter $\U$ is
 $\kappa$-uniform if for every set $U\in\U$ and any subset $K\subset U$ of size $|K|<\kappa$ the set
 $U\setminus K$ still belongs to $\U$.

To each $\kappa$-uniform ultrafilter $\U\in U_\kappa(A)$ assign the invariant filter $\F_\U=\bigcap_{x\in X}x\U$.
 This filter can be extended to a maximal invariant linked system $\mathcal L_\U$.
 We claim that $\mathcal L_\U\ne\mathcal L_\V$ for two different $\kappa$-uniform ultrafilters
 $\U,\V$ on $A$. Indeed, $\U\ne\V$ yields a subset $U\subset A$ such that $U\in\U$ and $U\notin\V$.
 Let $V=A\setminus U$. Since $\U$, $\V$ are $\kappa$-uniform, $|U|=|V|=\kappa$.

For every $\alpha<\kappa$ consider the sets $U_\alpha=\{a_\beta\in U:\beta>\alpha\}\in\U$ and
 $V_\alpha=\{a_\beta\in V:\beta>\alpha\}\in\V$.

It is clear that $$F_U=\bigcup_{\alpha<\kappa}x_\alpha U_\alpha\in\F_\U\mbox{ and }F_V=\bigcup_{\alpha<\kappa}x_\alpha V_\alpha\in\F_\V.$$

Let us show that $F_U\cap F_V=\emptyset$. Otherwise there would exist two ordinals $\alpha,\beta$
and points $u\in U_\alpha$, $v\in V_\beta$ such that $x_\alpha u= x_\beta v$. It follows from $u\ne v$ that $\alpha\ne\beta$.
Write the points $u,v$ as $u=a_\gamma$ and $v=a_\delta$ for some $\gamma> \alpha$ and $\delta>\beta$.
 Then we have the equality $x_\alpha a_\gamma=x_\beta a_\delta$. The inequality $u\ne v$implies that $\gamma\ne \delta$.
  We lose no generality assuming that
  $\delta>\gamma$. Then $$a_\delta=x_\beta^{-1}x_\alpha a_\gamma\in B_\delta^{-1}B_\delta A_{<\delta}$$
  which contradicts the choice of $a_\delta$.

Therefore, $F_U\cap F_V=\emptyset$. Taking into account that the
linked systems $\mathcal L_\U\supset \F_\U\ni F_U$ and $\mathcal
L_\V\supset \F_\V\ni F_V$ contain disjoint sets $F_U$, $F_V$, we
conclude that $\mathcal L_\U\ne\mathcal L_\V$. Consequently,
 $$|\inv[\lambda](X)|\ge|\{\mathcal L_{\U}:\U\in\U_\kappa(A)\}|=|U_\kappa(A)|=2^{2^{\kappa}}.$$
\end{proof}

The preceding theorem implies that $|\inv[\lambda](G)|=2^{\mathfrak c}$ for any countable group $G$.
Next, we evaluate the cardinality of $\inv[\lambda](G)$ for finite groups $G$.

Given a finite group $G$ consider the invariant linked system
$$\mathcal L_0=\{A\subset X:2|A|>|G|\}$$and the subset
$${\uparrow}\mathcal L_0=\{\A\in\inv[\lambda](G):\A\supset\mathcal L_0\}$$of $\inv[\lambda](G)$.

\begin{proposition}\label{p2.3} Let $G$ be a finite group. If $sl(G)\ge |G|/2$,
then $\inv[\lambda](G)={\uparrow}\mathcal L_0$.
\end{proposition}

\begin{proof} We should prove that each maximal invariant linked system $\A\in\inv[\lambda](G)$ contains $\mathcal L_0$.
Take any set $L\in\mathcal L_0$. Taking into account that $sl(G)\ge|G|/2$ and each set $A\in\A$ is self-linked, we
conclude that $|A|\ge|G|/2$ and hence $A$ intersects each shift $xL$ of $L$ (because $|A|+|xL|>|G|$). Since the set
 $L$ is self-linked, we get that the invariant linked system $\A\cup\{xL:x\in G\}$ is equal to $\A$ by the
 maximality of $\A$. Consequently, $L\in\A$ and hence $\mathcal L_0\subset\A$.
\end{proof}

In light of Proposition~\ref{p2.3} it is important to evaluate the cardinality of the set
${\uparrow}\mathcal L_0$. In $|G|$ is odd, then the invariant linked system $\mathcal L_0$ is maximal linked and
thus ${\uparrow}\mathcal L_0$ is a singleton.
The case of even $|G|$ is less trivial.

Given an group $G$ of finite even order $|G|$, consider the family
$$\mathcal S=\{A\subset G:AA^{-1}=G,\; |A|=|G|/2\}$$of self-linked subsets $A\subset G$
 of cardinality $|A|=|G|/2$. On the family $\mathcal S$ consider the equivalence
 relation $\sim$ letting $A\sim B$ for $A,B\in\mathcal S$ if there is
 $x\in G$ such that $A=xB$ or $X\setminus A=xB$. Let $\mathcal S/_\sim$ the quotient
 set of $\mathcal S$ by this equivalence relation and $s=|\mathcal S/_\sim|$ stand for the cardinality of $S/_\sim$.

\begin{proposition}\label{p2.4}  $|\inv[\lambda](G)|\ge|{\uparrow}\mathcal L_0|=2^s$.
\end{proposition}

\begin{proof} First we show that $\sim$ indeed is an equivalence relation on $\mathcal S$.
 So, assume that $\mathcal S\ne\emptyset$. Let us show that
$G\setminus A\in\mathcal S$ for every $A\in\mathcal S$. Let
$B=G\setminus A$. Assuming that $B\notin\mathcal S$, we conclude
that $B\cap x B=\emptyset$ for some $x\in G$. Since
$|B|=|A|=|G|/2$, we conclude that $xB=A$ and $G\setminus
A=B=x^{-1}A$. The equality $A\cap x^{-1}A=\emptyset$ implies
$x^{-1}\notin AA^{-1}=G$, which is a contradiction.

Taking into account that  $A=eA$ for every $A\in\mathcal S$, we conclude that $\sim$  is a
reflexive relation on $\mathcal S$. If $A\sim B$, then there is
 $x\in X$ such that $A=xB$ or $G\setminus A=xB$. This implies that
 $B=x^{-1}A$ or $X\setminus B=x^{-1}A$, that is $B\sim A$ and $\sim$ is
 symmetric. It remains to prove that the relation $\sim$ is transitive on $\mathcal S$. So let $A\sim B\sim C$.
  This means that there
 exist $x,y\in G$ such that $A=xB$ or $G\setminus A=xB$ and $B=yC$ or $G\setminus
 B=yC$. It is easy to check that in these cases $A=xyC$ or $X\setminus A=xyC$.

 Choose a subset $\mathcal T$ of $\mathcal S$ intersecting each equivalence class of $\sim$
  at a single point. Observe that $|\mathcal T|=|\mathcal S/_\sim|=s$. Now for every function
   $f:\mathcal T\to 2=\{0,1\}$ consider the maximal invariant linked system
$$\mathcal L_f=\mathcal L_0\cup\{xT:x\in G,\;T\in f^{-1}(0)\}\cup\{x(G\setminus T):x\in G,\;T\in f^{-1}(1)\}.$$
It can be shown that $$|{\uparrow}\mathcal L_0|=|\{\mathcal
L_f:f\in 2^{\mathcal T}\}|=2^{|\mathcal T|}=2^s.$$
\end{proof}

This proposition will help us to calculate the cardinality of the
set $\inv[\lambda](G)$ for all finite groups $G$ of order $|G|\le
8$:

\begin{theorem}\label{fincar2} The cardinality of $\inv[\lambda](G)$ for a group $G$ of size
$|G|\le 8$ can be found from the table:
{\small
$$
\begin{array}{|c|c|c|cc|c|cc|c|ccccc|}
\hline
G& C_2& C_3 &C_4 &C_2\oplus C_2& C_5 &D_6& C_6& C_7 & C_2^3 &D_8 &C_4\oplus C_2 & C_8 & Q_8\\
\hline
sl(G)&2&2&3&3&3&4&3&3&5&4&4&4&4\\
\inv[\lambda](X)&1&1&1&1&1&1&2&3&1&2&4&8&8\\
\hline
\end{array}
$$
}
\end{theorem}

\begin{proof} We divide the proof into $5$ cases.

1. If $sl(G)>|G|/2$, then $\mathcal L_0$ is a unique maximal invariant
 linked system and thus $|\inv[\lambda](X)|=1$. By Theorem~\ref{upbound}, $sl(G)>|G|/2$ if and
  only if $|G|\le 5$ or $G$ is isomorphic to $D_6$ or $C_2^3$.

2. If $sl(G)=|G|/2$, then $|\inv[\lambda](G)|=2^s$ where $s=|\mathcal S/_\sim|$. So it remains to
 calculate the number $s$ for the groups $C_6$, $D_8$, $C_4\oplus C_2$, $C_8$, and $Q_8$.

2a. If $G$ is cyclic of order 6, then we can take any generator $a$ on
$G$ and by routine calculations, check that
$$\mathcal S=\{x T,x(G\setminus T):x\in G\}$$ where $T=\{e,a,a^3\}$. It follows
that $s=|\mathcal S/_\sim|=1$ and thus $$|\inv[\lambda](G)|=|{\uparrow}\mathcal L_0|=2^s=2.$$

2b. If $G$ is cyclic of order 8, then we can take any generator
$a$ on $G$ and by routine verification check that
$$\mathcal S=\{xA,G\setminus xA,xB,G\setminus xB,C,G\setminus xC:x\in G\}$$where
$A=\{e,a,a^2,a^4\}$, $B=\{e,a,a^2,a^5\}$, and $C=\{e,a,a^3,a^5\}$. It follows that $s=|\mathcal S/_\sim|=3$ and
thus $$|\inv[\lambda](G)|=|{\uparrow}\mathcal L_0|=2^s=8.$$

2c. Assume that the group $G$ is isomorphic to $C_4\oplus C_2$ and
let $G_2=\{x\in G:xx=e\}$ be the Boolean subgroup of $G$. We claim
that a 4-element subset $A\subset G$ is self-linked if and only if
$|A\cap G_2|$ is odd.

To prove the ``if'' part of this claim, assume that $|A\cap
G_2|=3$. We claim that $A$ is self-linked. Let $A_2=A\cap G_2$ and
note that $G_2=A_2A^{-1}_2\subset AA^{-1}$ because
$|A_2|=3>2=|G_2|/2$. Now take any element $a\in A\setminus G_2$
and note that $AA^{-1}\supset aA_2^{-1}\cup A_2a^{-1}$. Observe
that both $aA_2^{-1}=aA_2$ and $A_2a^{-1}=a^{-1}A_2$ are 3-element
subsets in the 4-element coset $aG_2$. Those 3-element sets
 are different. Indeed, assuming that $aA_2^{-1}=A_2a^{-1}$ we would obtain that $a^2A_2=A_2$ which implies that
 $|A_2|=3$ is even. Consequently, $aG_2=aA_2^{-1}\cup A_2a^{-1}\subset AA^{-1}$ and finally $G=AA^{-1}$.

If $|A\cap G_2|=1$, then we can take any $a\in A\setminus G_2$ and consider the shift $Aa^{-1}$ which
has $|Aa^{-1}\cap G_2|=3$. Then the preceding case implies that $Aa^{-1}$ is self-linked and so is $A$.

To prove the ``only if'' part of the claim assume that $|A\cap G_2|$ is even. If $|A\cap G_2|=4$, then $A=G_2$
 and $AA^{-1}=G_2 G_2^{-1}=G_2\ne G$. If $|A\cap G_2|=0$, then $A=G_2a$ for any $a\in A$ and hence
 $AA^{-1}=G_2aa^{-1}G_2^{-1}=G_2\ne G$.
If $|A\cap G_2|=2$, then $|G_2\cap AA^{-1}|\le 3$ and again $AA^{-1}\ne G$.

Thus $$\mathcal S=\{A\subset G:\mbox{$|A|=4$ and $|A\cap G_2|$ is odd}\}.$$
Each set $A\in\mathcal S$ has a unique shift $aA$ with $aA\cap G_2=\{e\}$. There are exactly four subsets
$A\in\mathcal S$ with $A\cap G_2=\{e\}$ forming two equivalence classes with respect to the relation $\sim$.
Therefore $s=2$ and $$|\inv[\lambda](G)|=|{\uparrow}\mathcal L_0|=2^s=4.$$

2d. Assume that $G$ is isomorphic to the dihedral group $D_8$ of isometries of the square.
Then $G$ contains an element $a$ of order 4 generating a normal cyclic subgroup $H$. The element $a^2$
commutes with all the elements of the group $G$.

 We claim that for each self-linked 4-element subset $A\subset G$ we get $|A\cap H|=2$. Indeed, if $|A\cap H|$
 equals 0 or 4, then $A=Hb$ for some $b\in G$ and then $AA^{-1}=Abb^{-1}A^{-1}=H\ne G$. If $|A\cap H|$ equals 1 or 3,
 then replacing $A$ by a suitable shift, we can assume that $A\cap H=\{e\}$ and hence $A=\{e\}\cup B$ for
 some 3-element subset $B\subset G\setminus H$.
 It follows that $G\setminus H=AA^{-1}\setminus H=(B\cup B^{-1})=B\ne G\setminus H$. This contradiction shows
 that $|A\cap H|=2$. Without loss of generality, we can assume that $A\cap H=\{e,a^2\}$ (if it is not the case,
  replace $A$ by its shift $Ax^{-1}$ where $x,y\in A$ are such that $yx^{-1}=a^2$). Now take any element
  $b\in A\setminus H$. Since $G$ is not commutative, we get $ab=ba^3$. Observe that $ba^2\notin A$
  (otherwise $A=\{e,b,a^2,ba^2\}$ would be a subgroup of $G$ with $AA^{-1}=A\ne G$). Consequently,
  the 4-th element $c\in A\setminus\{e,a^2,b\}$ of $A$ should be of the form $c=ba$ or $c=ba^3=ab$. Observe
  that both the sets $A_1=\{e,a^2,b,ba\}$ and $A_2=\{e,a^2,b,ab\}$ are self-linked. Observe also that
  $$a^3(G\setminus A_1)=a^3\cdot\{a,a^3,ba^2,ba^3\}=\{e,a^2,ab,b\}=A_2.$$
Consequently, $s=|\mathcal S/_\sim|=1$ and $|\inv[\lambda](G)|=2^s=2$.

2e. Finally assume that $G$ is isomorphic to the group
$Q_8=\{\pm1,\pm i,\pm j,\pm k\}$ of quaternion units. The two-element subset $H=\{-1,1\}$ is a normal
subgroup in $X$. Let $\mathcal S_\pm=\{A\in S:H\subset A\}$ and observe that each set $A\in\mathcal S$
has a left shift in $\mathcal S$. Take any set $A\in\mathcal S_\pm$ and pick a point  $a\in A\setminus\{1,-1\}$.
Observe that the 4-th element $b\in A\setminus\{1,-1,a\}$ of $A$ is not equal to $-a$ (otherwise, $A$ is a
subgroup of $G$).

Conversely, one can easily check that each set $A=\{1,-1,a,b\}$
with $a,b\in G\setminus H$ and $a\ne-b$ is self-linked. This means
that $$\mathcal S_\pm=\{\{-1,1,a,b\}:a\ne-b\mbox{ and }a,b\in
G\setminus H\}$$and thus $|\mathcal S_\pm|=C^2_6-3=12$. Observe
that for each $A\in\mathcal S_2$ the set $-A\in\mathcal S_2$ and
there axactly two shifts of $X\setminus A$ that belong to
$\mathcal S_2$. This means that the equivalence class $[A]_\sim$
of any set $A\in\mathcal S$ intersects $\mathcal S_2$ in four
sets. Consequently, $s=|\mathcal S/_\sim|=|\mathcal
S_{\pm}|/4=12/4=3$ and
$$|\inv[\lambda](G)|=|{\uparrow}\mathcal L_0|=2^s=8.$$

3. If $|G|=7$, then $\mathcal L_0$ is one of three elements of
$\inv[\lambda](G)$. The other two elements can be found as
follows. Consider the invariant linked system
$$\mathcal L_1=\{A\subset G:|A|\ge5\}$$ and observe that $\mathcal L_1\subset\A$ for each $\A\in\inv[\lambda](G)$.
 Indeed, assuming that some $A\in\mathcal L_1$ does not belong to $\A$, we
 would conclude that $B=G\setminus A\in\A$ by the maximality of $\A$.
 Since $|G\setminus B|\le2$ we can find $x\in G\setminus BB^{-1}$. It
 follows that $B,xB$ are two disjoint sets in $\A$ which is not possible. Thus $\mathcal L_1\subset\A$.

Observe that $\mathcal L_1\subset\A\subset\mathcal L_0\cup\mathcal
L_3$, where
$$\mathcal L_3=\{A\subset G:|A|=3, \; AA^{-1}=G\}.$$ Given a
generator $a$ of the cyclic group $G$, consider the 3-element set
$T=\{a,a^2,a^4\}$ and note that $TT^{-1}=G$ and $T^{-1}\cap
T=\emptyset$. By a routine calculation, one  can check that
$$\mathcal L_3=\{xT,xT^{-1}:x\in G\}.$$ Since $T$ and $T^{-1}$ are
disjoint, the invariant linked system $\A$ cannot contain both the
sets $T$ and $T^{-1}$. If $\A$ contains none of the sets
$T,T^{-1}$, then $\A=\mathcal L_0$. If $\A$ contains $T$, then
$$\A=(\mathcal L_0\cup\{xT:x\in G\})\setminus\{y(G\setminus T):y\in
G\}.$$ If $T^{-1}\in \A$, then
$$\A=(\mathcal L_0\cup\{xT^{-1}:x\in G\})\setminus\{y(G\setminus T^{-1}):y\in G\}.$$
And those are the unique 3 maximal invariant systems in $\inv[\lambda](G)$.
\end{proof}

In the following theorem we characterize groups possessing a unique maximal invariant linked system.

\begin{theorem}\label{fingr} For a finite group $G$ the following conditions are equivalent:
\begin{enumerate}
\item $|\inv[\lambda](G)|=1$;
\item $sl(G)>|G|/2$;
\item $|G|\le 5$ or else $G$ is  isomorphic to $D_6$ or $C_2^3$.
\end{enumerate}
\end{theorem}

\begin{proof} $(2)\Ra(1)$. If $sl(G)>|G|/2$, then $\mathcal L_0=\{A\subset G:|A|>|G|/2\}$ is a unique
 maximal invariant linked system on $G$ (because invariant linked systems compose of self-linked sets).

$(1)\Ra(2)$ Assume that $sl(G)\le |G|/2$ and take a self-linked subset $A\subset G$ with $|A|\le|G|/2$.
 If $|G|$ is odd, then $\mathcal L_0$ is maximal linked and then any maximal invariant linked system
  $\A$ containing the self-linked set $A$ is distinct from $\mathcal L_0$, witnessing that $|\inv[\lambda](G)|>1$.

If $G$ is even, then we can enlarge $A$, if necessary, and assume that $|A|=|G|/2$.
We claim that the complement  $B=G\setminus A$ of $A$ is self-linked too. Assuming the converse, we would
 find some $x\notin BB^{-1}$ and conclude that $B\cap xB=\emptyset$, which implies that $A=G\setminus B=xB$
 and hence $x^{-1}A=B$. Then the sets $A$ and $x^{-1}A$ are disjoint which contradicts $x^{-1}\in AA^{-1}=G$.
 Thus $BB^{-1}=G$ which implies that $\{xB:x\in G\}$ is an invariant linked system. Since $|G|=2|A|$ is even,
 the unions $\A=\{xA:x\in G\}\cup\mathcal L_0$ and $\mathcal B=\{xB:x\in G\}\cup\mathcal L_0$ are invariant
 linked systems that can be enlarged to maximal linked systems $\tilde{\A}$ and $\tilde{\mathcal B}$, respectively.
Since the sets $A\in\mathcal A\subset\tilde{\A}$ and $B\in\mathcal B\subset\tilde{\mathcal B}$ are disjoint,
 $\tilde{\mathcal A}\ne\tilde{\mathcal B}$ are two distinct maximal invariant systems on $G$ and thus
 $|\inv[\lambda](G)|\ge 2$.

The equivalence $(2)\Leftrightarrow(3)$ follows from Theorem~\ref{upbound}(i).
\end{proof}

\section{Right zeros in $\lambda(X)$}\label{rzeros}

In this section we return to studying the superextensions of groups and shall detect groups $X$ whose superextensions
$\lambda(X)$ have right zeros. We shall show that for every group $X$ the right zeros of $\lambda(X)$
coincide with invariant maximal linked systems.

We recall that an element $z$ of a semigroup $S$ is called a {\em right} (resp.
{\em left}) {\em zero} in $S$ if $xz=z$ (resp. $zx=z$) for every
$x\in S$. This is equivalent to saying that the singleton $\{x\}$
is a left (resp. right) ideal of $S$.

By \cite[5.1]{G2} an inclusion hyperspace $\A\in G(X)$ is a right
zero in $G(X)$ if and only if $\A$ is invariant. This implies that
the minimal ideal of the semigroup $G(X)$ coincides with the set
$\invG(X)$ of invariant inclusion hyperspaces and is a compact
rectangular topological semigroup. We recall that a semigroup $S$
is called {\em rectangular} if $xy=y$ for all $x,y\in S$.

A similar characterization of right zeros holds also for the semigroup $\lambda(X)$.

\begin{proposition}\label{rightzero} A maximal linked system $\mathcal L$ is a right zero of
 the semigroup $\lambda(X)$ if and only if $\mathcal L$ is invariant.
\end{proposition}

\begin{proof} If $\mathcal L$ is invariant, then by proposition~5.1 of \cite{G2},
$\mathcal L$ is a right zero in $G(X)$ and consequently, a right zero in $\lambda(X)$.

Assume conversely that $\mathcal L$ is a right zero in $\lambda(X)$. Then for every
$x\in X$ we get $x\mathcal L=\mathcal L$, which means that $\mathcal L$ is invariant.
\end{proof}

Unlike the semigroup $G(X)$ which always contains right zeros, the
semigroup $\lambda(X)$ contains right zeros only for so-called odd
groups. We define a group $X$ to be {\em odd}
 if each element $x\in X$ has odd order. We recall that the {\em order} of an element $x$
 is the smallest integer number $n\ge 1$ such that $x^n$ coincides with the neutral element $e$ of $X$.

\begin{theorem}\label{odd} For a group $X$ the following conditions are equivalent:
\begin{enumerate}
\item the semigroup $\lambda(X)$ has a right zero;
\item some maximal invariant linked system on $X$ is maximal linked  (which can be written
as $\inv[\lambda](X)\cap\lambda(X)\ne\emptyset$);
\item each maximal invariant linked system is maximal linked  (which can be written as
 $\inv[\lambda](X)\subset\lambda(X)$);
\item for any partition $X=A\cup B$ either $AA^{-1}=X$ or $BB^{-1}=X$;
\item each element of $X$ has odd order.
\end{enumerate}
\end{theorem}

\begin{proof} The equivalence $(1)\Leftrightarrow(2)$ follows from Proposition~\ref{rightzero}.

$(2)\Ra(4)$ Assume that $\lambda(X)$ contains an invariant maximal linked system $\A$.
 Given any partition $X=A_1\cup A_2$, use the maximality of $\A$ to find $i\in\{1,2\}$ with $A_i\in\A$.
  We claim that $A_iA_i^{-1}=X$. Indeed, for every $x\in X$ the invariantness of $\A$ implies
   that $xA_i\in\A$ and hence $A_i\cap xA_i\ne\emptyset$, which implies $x\in A_iA_i^{-1}$.
\smallskip

$(4)\Ra(3)$ Assume that for every partition $X=A\cup B$ either
$AA^{-1}=X$ or $BB^{-1}=X$. We need to check that each maximal
invariant linked system $\mathcal L$ is maximal linked.
 In the other case, there would exist a set $A \in\mathcal L^\perp\setminus\mathcal L$.
 Since $\mathcal L\not\ni A$ is maximal invariant linked system, some shift $xA$ of $A$ does
  not intersect $A$ and thus $x\notin AA^{-1}$. Then our assumption implies that $B=X\setminus A$
  has property $BB^{-1}=X$, which means that the family $\{xB:x\in X\}$ is linked. We claim that
   $B\in\mathcal L^\perp$. Assuming the converse, we would find a set $L\in\mathcal L$ with
    $L\cap B=\emptyset$ and conclude that $A\in\mathcal L$ because $L\subset X\setminus B=A$.
    But this contradicts the choice of $A\in\mathcal L^\perp\setminus\mathcal L$.
     Therefore $B\in\mathcal L^\perp$ and $$\mathcal L\cup\{L\subset X:\exists x\in X\; (xB\subset L)\}$$
      is an invariant linked system that enlarges $\mathcal L$. Since $\mathcal L$ is a maximal invariant
      linked system, we conclude that $B\in\mathcal L$, which is not possible because $B$ does not
      intersect $A\in\mathcal L^\perp$. The obtained contradiction shows that
       $\mathcal L^\perp\setminus\mathcal L=\emptyset$, which means that $\mathcal L$ belongs to
       $\lambda(X)$ and thus is an invariant maximal linked system.
\smallskip

The implication $(3)\Ra(2)$ is trivial.
\smallskip

$\neg(5)\Ra\neg(4)$ Assume that $X\setminus\{e\}$ contains a point $a$ whose order is even or infinity.
 Then the cyclic subgroup $H=\{a^n:n\in\IZ\}$ generated by $a$ decomposes into two disjoint sets
 $H_1=\{a^n:n\in2\IZ+1\}$ and $H_2=\{a^n:n\in2\IZ\}$
such that $aH_1=H_2$. Take a subset $S\subset X$ meeting each
coset $Hx$, $x\in X$, in a single point. Consider the disjoint
sets $A_1=H_1S$ and $A_2=H_2S$ and note that
$aA_1=A_2=X\setminus A_1$ and $aA_2=X\setminus A_2$, which implies that $a\notin A_iA_i^{-1}$ for $i\in\{1,2\}$.
 Since $A_1\cup A_2=X$, we get a negation of (4).
\smallskip

$(5)\Ra(4)$ Assume that each element of $X$ has odd order and assume that $X$ admits a partition
 $X=A\sqcup B$ such that $a\notin AA^{-1}$ and $b\notin BB^{-1}$ for some $a,b\in X$. Then $aA\subset X\setminus A=B$
  and $bB\subset X\setminus B=A$. Observe that
$$baA\subset bB\subset A$$and by induction, $(ba)^iA\subset A$ for all $i>0$. Since all elements
 of $X$ have finite order, $(ba)^n=e$ for some $n\in\IN$. Then $(ba)^{n-1}A\subset A$ implies
$$A=(ba)^nA\subset baA\subset bB\subset A$$ and hence $bB=A$.
It follows from $$X=bA\sqcup bB=bA\sqcup A=B\sqcup A$$ that $bA=B$. Thus $x\in A$ if and only if $bx\in B$.

 Let $H=\{b^n:n\in\IZ\}\subset X$ be the cyclic subgroup generated by $b$. By our assumption it
 is of odd order. On the other hand, the equality $bB=A=b^{-1}B$ implies that the intersections $H\cap A$
 and $H\cap B$ have the same cardinality because $b(B\cap H)=A\cap H$. But this is not possible
 because of the odd cardinality of $H$.
\end{proof}

\section{(Left) zeros of the semigroup $\lambda(X)$}\label{lzeros}

An element $z$ of a semigroup $S$ is called a {\em zero} in $S$ if $xz=z=zx$ for all $x\in S$.
This is equivalent to saying that $z$ is both a left and right zero in $S$.

\begin{proposition}\label{zero} Let $X$ be a group. For a maximal linked system $\mathcal L\in\lambda(X)$
 the following conditions are equivalent:
\begin{enumerate}
\item $\mathcal L$ is a left zero in $\lambda(X)$;
\item $\mathcal L$ is a zero in $\lambda(X)$;
\item $\mathcal L$ is a unique invariant maximal linked system on $X$.
\end{enumerate}
\end{proposition}

\begin{proof} $(1)\Ra(3)$ Assume that $\mathcal Z$ is a left zero in $\lambda(X)$.
Then $\mathcal Zx=\mathcal Z$ for all $x\in X$ and thus $$\mathcal
Z^{-1}=\{Z^{-1}:Z\in\mathcal Z\}$$ is an invariant maximal linked
system on $X$, which implies that the group $X$ is odd according to Theorem~\ref{odd}. Note that
 for every right zero $\mathcal A$ of $\lambda(X)$ we get
$$\mathcal Z=\mathcal Z\circ\mathcal A=\mathcal A$$ which implies that $\mathcal Z$ is a unique right
 zero in $\lambda(X)$ and by Proposition~\ref{rightzero} a unique invariant maximal linked system on $X$.
\smallskip

$(3)\Ra(2)$ Assume that $\mathcal Z$ is a unique invariant maximal linked system on $X$.
We claim that $\mathcal Z$ is a left zero of $\lambda(X)$. Indeed, for every $\A\in\mathcal A$
 and $x\in X$ we get $x\mathcal Z\circ\mathcal A=\mathcal Z\circ\mathcal A$, which means that $\Z\circ\A$
  is an invariant maximal linked system. By Proposition \ref{rightzero}, $\mathcal Z\circ\A$ is a right zero and hence
  $\mathcal Z\circ\A=\mathcal Z$ because $\mathcal Z$ is a unique right zero. This means that $\mathcal Z$ is
   a left zero, and being a right zero, a zero in $\lambda(X)$.

$(2)\Ra(1)$ is trivial.
\end{proof}

\begin{theorem}\label{zerochar} The superextension $\lambda(X)$ of a group $X$ has a zero if and only if $X$ is
 isomorphic to $C_1$, $C_3$ or $C_5$.
\end{theorem}

\begin{proof} If $X$ is a group of odd order $|X|\le 5$, then $\inv[\lambda](X)\subset\lambda(X)$ because
$X$ is odd and $|\inv[\lambda](X)|=1$ by Theorem~\ref{fingr}.
This means that $\lambda(X)$ contains a unique invariant maximal
linked system, which is the zero of $\lambda(X)$ by
Proposition~\ref{zero}.

Now assume conversely that the semigroup $\lambda(X)$ has a zero
element $\mathcal Z$. By Proposition~\ref{rightzero} and
Theorem~\ref{odd}, $X$ is odd and thus
$\inv[\lambda](X)\subset\lambda(X)$. Since the
 zero $\mathcal Z$ of $\lambda(X)$ is a unique invariant maximal linked system on $X$,
 we get $|\inv[\lambda](X)|\le 1$. By Theorem~\ref{fingr}, $X$ has order $|X|\le 5$ or is
 isomorphic to $D_3$ or $C_2^3$. Since $X$ is odd, $X$ must be isomorphic to $C_1$, $C_3$ or $C_5$.
\end{proof}

\section{The commutativity of $\lambda(X)$}\label{scom}

In this section we detect groups $X$ with commutative superextension.

\begin{theorem}\label{comchar} The superextension $\lambda(X)$ of a group $X$ is commutative if and only if $|X|\le 4$.
\end{theorem}

\begin{proof} The commutativity of the superextensions $\lambda(X)$ of groups $X$ of order $|X|\le 4$ will
 be established in Section~\ref{fin}.

Now assume that a group $X$ has commutative superextension $\lambda(X)$.
Then $X$ is commutative. We need to show that $|X|\le 4$. First we show that $|\inv[\lambda](X)|=1$.

Assume that $\inv[\lambda](X)$ contains two distinct maximal invariant linked systems $\A$ and $\mathcal B$.
 Taking into account that $\A,\mathcal B\in\inv[\lambda](X)\subset\inv[G](X)$ and each element of $\inv[G](X)$
 is a right zero in $G(X)$ (see  \cite[5.1]{G2}) we conclude that
 $$\A\circ \mathcal B=\mathcal B\ne\A=\mathcal B\circ\A.$$
Extend the linked systems systems $\A,\mathcal B$ to maximal linked systems $\tilde{\A}\supset\A$
 and $\tilde{\mathcal B}\supset\mathcal B$. Because of the commutativity of $\lambda(X)$, we get
$$\A=\mathcal B\circ\mathcal A\subset\tilde{\mathcal B}\circ\tilde{\mathcal A}=\tilde{\mathcal A}\circ\tilde{\mathcal B}\supset\A\circ\mathcal B=\mathcal B.$$
This implies that the union $\A\cup\mathcal B\ne\A$ is an invariant linked system extending $\A$,
which is not possible because of the maximality of $\A$. This contradiction shows that $|\inv[\lambda](X)|=1$.
 Applying Theorem~\ref{fingr}, we conclude that $|X|\le 5$ or $X$ is isomorphic to $C_2^3$.

It remains to show that the semigroups $\lambda(C_5)$ and
$\lambda(C_2^3)$ are not commutative. The non-commutativity of
$\lambda(C_5)$ will be shown in Section~\ref{fin}.

To see that $\lambda(C_2^3)$ is not commutative, take any 3
generators $a,b,c$ of $C_2^3$ and consider the sets
$A=\{e,a,b,abc\}$, $H_1=\{e,a,b,ab\}$, $H_2=\{e,a,bc,abc\}$.
Observe that $H_1,H_2$ are subgroups in $C_2^3$. For every
$i\in\{1,2\}$ consider the linked system
$\A_i=\langle\{H_1,H_2\}\cup\{xA:x\in H_i\}\rangle$ and extend it
to a maximal linked system $\tilde\A_i$ on $C_2^3$.

We claim that the maximal linked systems $\tilde\A_1$ and $\tilde
\A_2$ do not commute. Indeed,
$$\tilde\A_2\circ\tilde\A_1\ni\bigcup_{x\in
H_1}x*(x^{-1}bA)=bA=\{e,b,ba,ac\},$$
$$\tilde\A_1\circ\tilde\A_2\ni\bigcup_{x\in
H_2}x*(x^{-1}bcA)=bcA=\{a,c,bc,abc\}.$$ It follows from $bA\cap
bcA=\emptyset$ that
$\tilde\A_1\circ\tilde\A_2\neq\tilde\A_2\circ\tilde \A_1$.
\end{proof}

\section{The superextensions of finite groups}\label{fin}

In this section we shall describe the structure of the superextensions $\lambda(G)$  of finite groups
 $G$ of small cardinality (more precisely, of cardinality $|G|\le5$).  It is known that the cardinality
of $\lambda(G)$ growth very quickly as $|G|$ tends to infinity.
The calculation of the cardinality of $|\lambda(G)|$ seems to be a
difficult combinatorial problem related to the still unsolved
Dedekind's problem of calculation of the number $M(n)$ of
inclusion hyperpspaces on an $n$-element subset, see \cite{De}. We
were able to calculate the cardinalities of $\lambda(G)$ only
  for groups $G$ of cardinality $|G|\le6$. The results of (computer) calculations
  are presented in the following table:
$$
\begin{array}{r|cccccc}
|G| &   1 & 2 & 3 & 4 & 5 & 6\\
\hline
|\lambda(G)|& 1 & 2 & 4 & 12&81 &2646\\
|\lambda(G)/G| &1 &1 &2 & 3 & 17 & 447
\end{array}
$$

Before describing the structure of superextensions of finite groups, let us make some remarks
concerning the structure of a semigroup $S$ containing a group $G$. In this case
$S$ can be thought as a $G$-space endowed with the left action
of the group $G$. So we can consider the orbit space
$S/G=\{Gs:s\in S\}$ and the projection $\pi:S\to S/G$. If $G$ lies
in the center of the semigroup $S$ (which means that the elements
of $G$ commute with all the elements of $S$), then the orbit space
$S/G$ admits a unique  semigroup operation turning $S/G$ into a
semigroup and the orbit projection $\pi:S\to S/G$ into a semigroup
homomorphism. A subsemigroup $T\subset S$ will be called a {\em
transversal semigroup} if the restriction $\pi:T\to S/G$ is an
isomorphism of the semigroups. If $S$ admits a transversal
semigroup $T$, then it is a homomoprhic image of the product
$G\times T$ under the semigroup homomorphism
$$h:G\times T\to S,\;\; h:(g,t)\mapsto gt.$$
This helps to recover the algebraic structure of $S$ from the
structure of a transversal semigroup.

For a system $\mathcal B$ of subsets of a set $X$ by $$\langle
\mathcal B\rangle=\{A\subset X:\exists B\in\mathcal B\;\;(B\subset
A)\}$$ we denote the inclusion hyperspace generated by $\mathcal
B$.
\smallskip

Now we shall analyse the entries of the above table. First note that
each group $G$ of size $|G|\le 5$ is abelian and is isomorphic to one of the groups: $C_1$,
$C_2$, $C_3$, $C_4$, $C_2\oplus C_2$, $C_5$.
It will be convenient to think of the cyclic group $C_n$ as the
multiplicative subgroups $\{z\in\IC:z^n=1\}$ of the complex plane.

\subsection{The semigroups $\lambda(C_1)$ and $\lambda(C_2)$} For the groups $C_n$ with
$n\in\{1,2\}$ the semigroup $\lambda(C_n)$ coincides with $C_n$ while the orbit semigroup
$\lambda(C_n)/C_n$ is trivial.

\subsection{The semigroup $\lambda(C_3)$}
For the group $C_3$ the semigroup $\lambda(C_3)$ contains the three principal ultrafilters $1,z,-z$
where $z=e^{2\pi i/3}$ and the maximal linked inclusion hyperspace
$\triangleright=\langle\{1,z\},\{1,-z\},\{z,-z\}\rangle$ which is the zero in $\lambda(C_3)$.
The superextension $\lambda(C_3)$ is isomorphic to the multiplicative semigroup  $C_3^0=\{z\in\IC:z^4=z\}$
of the complex plane.
The  latter semigroup has zero $0$ and unit $1$ which are the unique idempotents.

The transversal semigroup $\lambda(C_3)/C_3$ is isomorphic to the semilattice $2=\{0,1\}$ endowed with
the $\min$-operation.

 \subsection{The semigroups $\lambda(C_4)$ and $\lambda(C_2\oplus C_2)$} The semigroup $\lambda(C_4)$
  contains 12 elements while the orbit semigroup $\lambda(C_4)/C_4$ contains 3 elements. The semigroup
  $\lambda(C_4)$ contains a transversal semigroup $$\lambda_T(G)=\{1,\triangle,\square\}$$ where $1$
   is the neutral element of $C_4=\{1,-1,i,-i\}$,
$$
\begin{aligned}
\triangle=\;&\langle\{1,i\},\{1,-i\},\{i,-i\}\rangle\mbox{ and }\\
\square=\;&\langle\{1,i\},\{1,-i\},\{1,-1\},\{i,-i,-1\}\rangle.
\end{aligned}$$
The transversal semigroup is isomorphic to the extension $C_2^1=C_2\cup\{e\}$ of the cyclic group $C_2$ by an
external unit $e\notin C_2$ (such that $ex=x=xe$ for all $x\in C_2^1$). The action of the group $C_4$ on
$\lambda(C_4)$ is free so, $\lambda(C_4)$ is
isomorphic to $\lambda_T(C_4)\oplus C_4$.

The semigroup $\lambda(C_2\oplus C_2)$ has a similar algebraic structure. It contains a transversal
 semigroup $$\lambda_T(C_2\oplus C_2)=\{e,\triangle,\square\}\subset\lambda(C_2\oplus C_2)$$
  where $e$ is the principal ultrafilter supported by the neutral element $(1,1)$ of $C_2\oplus C_2$ and
the maximal linked inclusion hyperspaces $\triangle$ and $\square$ are defined by analogy with the
 case of the group $C_4$:
{\small
$$
\begin{aligned}
\triangle=\;&\langle\{(1,1),(1,-1)\},\{(1,1),(-1,1)\},\{(1,-1),(-1,1)\}\rangle\mbox{ and }\\
\square=\;&\langle\{(1,1),(1,-1)\},\{(1,1),(-1,1)\},\{(1,1),(-1,-1)\},\{(1,-1),(-1,1),(-1,-1)\}\rangle.
\end{aligned}$$
}
The transversal semigroup $\lambda_T(C_2\oplus C_2)$ is isomorphic to $C_2^1$ and $\lambda(C_2\oplus C_2)$
is isomorphic to $C_2^1\oplus C_2\oplus C_2$.

We summarize the obtained results on the algebraic structure of the semigroups $\lambda(C_4)$ and
$\lambda(C_2\oplus C_2)$ in the following proposition.

\begin{proposition} Let $G$ be a group of cardinality $|G|=4$.
\begin{enumerate}
\item The semigroup $\lambda(G)$ is isomorphic to $C_2^1\oplus
G$ and thus is commutative; \item $\lambda(G)$ contains two
idempotents; \item $\lambda(G)$ has a unique proper ideal
$\lambda(G)\setminus G$ isomorphic to the group $C_2\oplus G$.
\end{enumerate}
\end{proposition}

\subsection{The semigroup $\lambda(C_5)$.} Unlike to $\lambda(C_4)$, the semigroup
$\lambda(C_5)$ has  complicated algebraic structure. It contains 81 elements.
One of them is zero
$$\Z=\{A\subset C_5:|A|\ge 3\},$$ which is invariant under any bijection of $C_5$.
All the other 80 elements have 5-element orbits under the action of $C_5$,
which implies that the orbit semigroup $\lambda(C_5)/C_5$ consists of 17 elements.
Let $\pi:\lambda(C_5)\to\lambda(C_5)/C_5$ denote the orbit projection.

It will be convenient to think of $C_5$ as the field $\{0,1,2,3,4\}$ with the
multiplicative subgroup $C_5^*=\{1,-1,2,-2\}$ of invertible elements (here $-1$
and $-2$ are identified with $4$ and $3$, respectively). Also for elements $x,y,z\in C_5$ we shall
write $xyz$ instead of $\{x,y,z\}$.

The semigroup $\lambda(C_5)$ contains 5 idempotents:
$$\begin{aligned}
\U=&\langle 0\rangle,\;\Z,\\
\Lambda_4=&\langle 01,02,03,04,1234\rangle,\\
{\Lambda}=&\langle 02,03,123,014,234\rangle,\\
2{\Lambda}=&\langle 04,01,124,023,143\rangle,
\end{aligned}$$
which commute and thus form an abelian subsemigroup $E(\lambda(C_5))$.
Being a semilattice, $E(\lambda(C_5))$ carries a natural partial order:
$e\le f$ iff $e\circ f=e$. The partial order $$\Z\le {\Lambda},2{\Lambda}\le \Lambda_4\le\U$$
 on the set $E(\lambda(C_5))$ is designed at the picture:

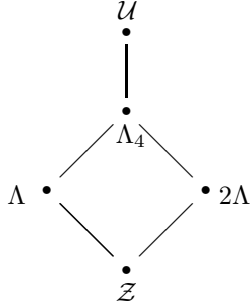
\begin{picture}(100,150)(-150,0)
\put(-4,18){$\Z$}
\put(0,30){\circle*{3}}
\put(5,35){\line(1,1){20}}
\put(-5,35){\line(-1,1){20}}
\put(-30,60){\circle*{3}}
\put(-45,55){${\Lambda}$}

\put(30,60){\circle*{3}}
\put(35,55){$2{\Lambda}$}
\put(-4,78){$\Lambda_4$}
\put(0,90){\circle*{3}}
\put(0,95){\line(0,1){20}}
\put(-5,35){\line(-1,1){20}}
\put(5,85){\line(1,-1){20}}
\put(-5,85){\line(-1,-1){20}}
\put(0,120){\circle*{3}}
\put(-3,125){$\U$}
\end{picture}

The other distinguished subset of $\lambda(C_5)$ is
$$
\begin{aligned}
\sqrt{E(\lambda(C_5))}=&\,\{\mathcal L\in\lambda(C_5):\mathcal L\circ\mathcal L\in E(\lambda(C_5))\}=\\
=&\,\{\mathcal L\in\lambda(C_5):\mathcal L\circ\mathcal L\circ\mathcal L\circ\mathcal L=\mathcal L\circ\mathcal L\}.
\end{aligned}
$$ We shall show that this set contains a point from each $C_5$-orbit in $\lambda(C_5)$.

First we show that this set has at most one-point intersection with each orbit.
 Indeed, if $\mathcal L\in\sqrt{E(\lambda(C_5))}$ and $\mathcal L\circ\mathcal L\ne\Z$,
 then  for every $a\in C_5\setminus\{0\}$, we get
$$
\begin{aligned}
(\mathcal L+a)\circ(\mathcal L+a)\circ(\mathcal L+a)\circ(\mathcal L+a)=&\,\mathcal L\circ\mathcal L\circ\mathcal L\circ\mathcal L+4a=\\
=&\,\mathcal L\circ\mathcal L+4a\ne \mathcal L\circ\mathcal L+2a=(\mathcal L+a)\circ(\mathcal L+a).
\end{aligned}
$$witnessing that $\mathcal L+a\notin\sqrt{\lambda_T(C_5)}$.

By a direct calculation one can check that the set $\lambda_T(C_5)$  contains the
following four maximal linked systems:$$\begin{aligned}
\Delta=&\langle 02,03,23\rangle,\\
\Lambda_3=&\langle 02,03,04,234\rangle,\\
\Theta=&\langle 14,012,013,123,024,034,234\rangle,\\
\Gamma=&\langle 02,04,013,124,234\rangle.
\end{aligned}
$$
For those systems we get
$$
\begin{aligned}
&\Delta\circ \Delta=\Delta\circ\Delta\circ\Delta={\Lambda},\\
&\Lambda_3\circ\Lambda_3=\Lambda_3\circ\Lambda_3\circ\Lambda_3={\Lambda},\\
&\F\circ\Theta=\F\circ\Gamma=\Z\;\; \mbox{for every $\F\in\lambda(C_5)\setminus C_5.$}
\end{aligned}
$$
All the other elements of $\lambda(C_5)$ can be found as  images of
$\Delta,\Theta,\Gamma,\Lambda_3$ under the affine transformations of the
 field $C_5$. Those are maps of the form $$f_{a,b}:x\mapsto ax+b \mod 5,$$
 where $a\in\{1,-1,2,-2\}=C_5^*$ and $b\in C_5$. The image of a maximal linked
  system $\mathcal L\in\lambda(C_5)$ under such a transformation will be denoted by $a\mathcal L+b$.

One can check that $a\Lambda_4=\Lambda_4$ for each $a\in C_5^*$
while ${\Lambda}=-{\Lambda}$, and $\Theta=-\Theta$. Since the
linear transformations of the form $f_{a,0}:C_5\to C_5$, $a\in
C_5^*$, are authomorphisms of the group $C_5$ the induced
transformations $\lambda f_{a,0}:\lambda(C_5)\to\lambda(C_5)$
are authomorphisms of the semigroup $\lambda(C_5)$. This implies
that those transformations do not move the subsets
$E(\lambda(C_5))$ and $\sqrt{E(\lambda(C_5))}$. Consequently,
the set $\sqrt{E(\lambda(C_5)}$ contains the maximal linked
systems:
$$a\Delta,a\Theta,a\Lambda_3,a\Gamma, \;a\in\IZ^*_5,$$ which together with the idempotents form a 17-element subset
$$T_{17}=E(\lambda(C_5))\cup\big\{a\Delta,a\Theta:a\in\{1,2\}\big\}\cup \{a\Lambda_3,a\Gamma:a\in\IZ^*_5\}$$
that projects bijectively onto the orbit semigroup $\lambda(C_5)/C_5$. The set $T_{17}$ looks as follows
 (we connect an element $x\in T_{17}$ with an idempotent $e\in T_{17}$ by an arrow if $x\circ x=e$):

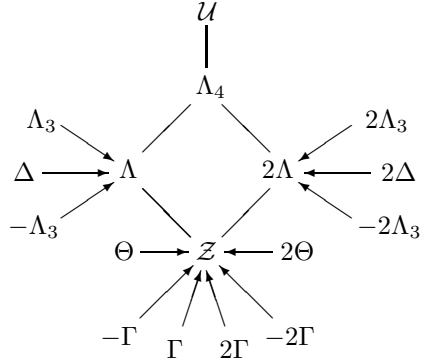
\begin{picture}(100,150)(-150,-20)
\put(-4,27){$\Z$}
\put(6,36){\line(1,1){18}}
\put(-6,36){\line(-1,1){18}}
\put(-33,58){$\Lambda$}

\put(21,58){$2\Lambda$}
\put(-4,90){$\Lambda_4$}
\put(0,100){\line(0,1){16}}
\put(-6,36){\line(-1,1){18}}
\put(6,86){\line(1,-1){18}}
\put(-6,86){\line(-1,-1){18}}
\put(-3,118){$\U$}

\put(-40,-5){$-\Gamma$}
\put(-35,27){$\Theta$}
\put(-25,30){\vector(1,0){18}}
\put(28,27){$2\Theta$}

\put(25,30){\vector(-1,0){18}}
\put(-25,5){\vector(1,1){20}}
\put(-15,-10){$\Gamma$}
\put(-9,2){\vector(1,3){7}}
\put(5,-10){$2\Gamma$}
\put(7,2){\vector(-1,3){7}}
\put(22,-5){$-2\Gamma$}
\put(25,5){\vector(-1,1){20}}

\put(-75,37){$-\Lambda_3$}
\put(-55,43){\vector(3,2){20}}
\put(-73,57){$\Delta$}
\put(-62,60){\vector(1,0){25}}
\put(-68,77){$\Lambda_3$}
\put(-55,78){\vector(3,-2){20}}

\put(57,37){$-2\Lambda_3$}
\put(55,43){\vector(-3,2){20}}
\put(66,57){$2\Delta$}
\put(62,60){\vector(-1,0){25}}
\put(60,77){$2\Lambda_3$}
\put(55,78){\vector(-3,-2){20}}
\end{picture}

The set $\sqrt{E(\lambda(C_5))}$ includes 24 elements more and coincides with the union
 $T_{17}\cup\sqrt{\Z}$ where $$\sqrt{\Z}=\{a\Theta+b,a\Gamma+b:a\in\IZ^*_5,\; b\in C_5\}.$$

Since each element of $\lambda(C_5)$ can be uniquely written as the sum $\mathcal L+b$ for
 some $\mathcal L\in T_{17}$ and  $b\in C_5$, the multiplication table for the semigroup
  $\lambda(C_5)$ can be recovered from the Cayley table for multiplication of the elements from $T_{17}$:
\medskip

{\small
$$
\begin{array}{|c|c|cccc|cccc|c|}

\hline
\circ & \Lambda_4 & {\Lambda} & \Delta & \Lambda_3 & -\Lambda_3 & 2{\Lambda} & 2\Delta & 2\Lambda_3 & -2\Lambda_3 & a\Theta,a\Gamma\\
\hline

\Lambda_4 & \Lambda_4 & {\Lambda} & {\Lambda}  & {\Lambda} & {\Lambda}  & 2{\Lambda} & 2{\Lambda} & 2{\Lambda} & 2{\Lambda} & \Z\\
\hline

{\Lambda} & {\Lambda} & {\Lambda} & {\Lambda}  & {\Lambda} & {\Lambda}  & \Z & \Z & \Z & \Z & \Z\\
\hline
\Delta & \Delta & {\Lambda} & {\Lambda}  & {\Lambda} & {\Lambda}  &
2\Theta & 2\Theta & 2\Theta & 2\Theta & \Z\\

\Lambda_3 & \Lambda_3 & {\Lambda} & {\Lambda}  & {\Lambda} & {\Lambda}  &
2\Theta+2 & 2\Theta+2 & 2\Theta+2 & 2\Theta+2 & \Z\\

-\Lambda_3 & -\Lambda_3 & {\Lambda} & {\Lambda} & {\Lambda} & {\Lambda}  &
2\Theta-2 & 2\Theta-2 & 2\Theta-2 & 2\Theta-2 & \Z\\
\hline

2{\Lambda} & 2{\Lambda} & \Z & \Z & \Z & \Z &
2{\Lambda} & 2{\Lambda}  & 2{\Lambda} & 2{\Lambda} & \Z\\
\hline

2\Delta & 2\Delta & \Theta & \Theta & \Theta & \Theta &
2{\Lambda} & 2{\Lambda}  & 2{\Lambda} & 2{\Lambda} & \Z\\

2\Lambda_3 & 2\Lambda_3 & \Theta-1 & \Theta-1 & \Theta-1 &
\Theta-1 &
2{\Lambda} & 2{\Lambda}  & 2{\Lambda} & 2{\Lambda} & \Z\\

-2\Lambda_3 & -2\Lambda_3 & \Theta+1 & \Theta+1 & \Theta+1 &
\Theta+1 &
2{\Lambda} & 2{\Lambda}  & 2{\Lambda} & 2{\Lambda} & \Z\\
\hline

\Theta & \Theta & \Theta & \Theta & \Theta & \Theta &
\Z & \Z & \Z & \Z & \Z\\

2\Theta & 2\Theta & \Z & \Z & \Z & \Z &
2\Theta & 2\Theta & 2\Theta & 2\Theta & \Z\\
\hline \Gamma & \Gamma & \Theta+1 & \Theta+1 & \Theta+1 & \Theta+1
&
2\Theta+2 & 2\Theta+2 & 2\Theta+2 & 2\Theta+2 & \Z\\

-\Gamma & -\Gamma & \Theta-1 & \Theta-1 & \Theta-1 & \Theta-1 &
2\Theta-2 & 2\Theta-2 & 2\Theta-2 & 2\Theta-2 & \Z\\

2\Gamma & 2\Gamma & \Theta-1 & \Theta-1 & \Theta-1 & \Theta-1 &
2\Theta+2 & 2\Theta+2 & 2\Theta+2 & 2\Theta+2 & \Z\\

-2\Gamma & -2\Gamma & \Theta+1 & \Theta+1 & \Theta+1 & \Theta+1 &
2\Theta-2 & 2\Theta-2 & 2\Theta-2 & 2\Theta-2 & \Z\\
\hline
\end{array}
$$
}

Looking at this table we can see that $T_{17}$ is not a
subsemigroup of $\lambda(C_5)$ and hence is not a transversal
semigroup for $\lambda(C_5)$. This is not occasional.

\begin{proposition} The semigroup $\lambda(C_5)$ contains no transversal semigroup.
\end{proposition}

\begin{proof} Assume conversely that $\lambda(C_5)$ contains a subsemigroup $T$
that projects bijectively onto the orbit semigroup
$\lambda(C_5)/C_5$. Then $T$ must include the set
$E(\lambda(C_5))$ of idempotents and also the subset
$\sqrt{E(\lambda(C_5))}\setminus\sqrt{\Z}$. Consequently,
$$T\supset\{\U,\Z,{\Lambda},-{\Lambda},\Delta,2\Delta,\Lambda_3,-\Lambda_3,2\Lambda_3,-2\Lambda_3\}.$$
Since
$2\Lambda_3\circ\Lambda=\Theta-1\neq\Theta=2\Delta\circ\Lambda$,
then there are two different points in the intersection
$T\cap(\Theta+C_5)$ which should be a singleton. This
contradiction completes the proof.

\end{proof}

Analysing the Cayley table for the set $T_{17}$ we can establish the following
properties of the semigroup $\lambda(C_5)$.

\begin{proposition}
\begin{enumerate}
\item The maximal linked system $\Z$ is the zero of $\lambda(\IZ)$.
\item $\lambda(C_5)$ contains 5 idempotents: $\U$, $\Z$, $\Lambda_4$, $\Lambda$, $2\Lambda$, which commute.
\item The set of central elements of $\lambda(C_5)$ coincides with $C_5\cup\{\Z\}$.
\item All non-trivial subgroups of $\lambda(C_5)$ are isomorphic to $C_5$.
\end{enumerate}
\end{proposition}

\subsection{Summary table} The obtained results on the superextensions of groups $G$ with $|G|\le 5$
are summed up in the following table in which $K(\lambda(G))$ stands for the minimal ideal of $\lambda(G)$.

$$
\begin{array}{|c|ccccc|}
\hline
|G| & |\lambda(G)| & \lambda(G) & |E(\lambda(G))| & K(\lambda(G)) & \mbox{maximal group}\\
\hline
2 & 2 & C_2 & 1 & C_2 &  C_2 \\
3 & 4 & C_3\cup\{\triangleright\} & 2 & \{\triangleright\} & C_3\\
4 & 12 & C_2^1\times G & 2 & C_2\times G & C_2\times G\\
5 & 81 & T_{17}\cdot C_5 & 5 & \{\Z\} & C_5 \\
\hline
\end{array}
$$


\end{document}